\documentclass[graybox]{svmult}

\usepackage{type1cm}        
\usepackage{makeidx}         
\usepackage{graphicx}        
\usepackage{multicol}        
\usepackage[bottom]{footmisc}

\usepackage{newtxtext}       %
\usepackage[varvw]{newtxmath}       

\makeindex             
\usepackage{url,doi,hyperref}

\usepackage{amsmath, amsthm, amsfonts}

\usepackage{hyperref}
\usepackage{epic}
\usepackage[dvips]{epsfig}
\usepackage{enumitem}
\usepackage{csquotes}

\definecolor{darkcyan}{rgb}{0.,0.5,0.5}
\definecolor{darkgreen}{rgb}{0.,0.7,0.}
\definecolor{brown}{rgb}{0.5,0.5,0.}

\theoremstyle{plain}

\theoremstyle{definition}

\setlist{itemsep=5pt}
\usepackage{comment}
\newcommand{\nLtwo}[1]{\|{#1}\|_{L^2(\Omega)}}
\newcommand{\prodLtwo}[2]{\langle{#1},{#2}\rangle_{L^2(\Omega)}}

\newcommand{\SMGT}{SMGT}

\newcommand{\revision}[1]{{ #1}}
\newcommand{\revisionown}[1]{{ #1}}
\newcommand{\MarginRone}[1]{}
\newcommand{\MarginRtwo}[1]{}

\begin{document}

\title*{On the Jordan-Moore-Gibson-Thompson equation of nonlinear acoustics}
\titlerunning{The Jordan-Moore-Gibson-Thompson equation}
\author{Barbara Kaltenbacher
}
\authorrunning{Kaltenbacher
}
\institute{Barbara Kaltenbacher \at 
Department of Mathematics,
Alpen-Adria-Universit\"at Klagenfurt, 
Universit\"atsstra\ss e 65-67, 
9020 Klagenfurt, Austria, \email{barbara.kaltenbacher@aau.at}
}
%
%
\maketitle

\abstract{The JMGT equation was put forward by Pedro Jordan~\cite{jordan2008nonlinear,jordan2014second}, also referring to earlier work by Moore and Gibson~\cite{moore1960propagation}, as well as Thompson~\cite{thompson} to amend the infinite speed of sound paradox of classical models of nonlinear acoustics such as the 
\revision{Westervelt and Kuznetsov equations.}
\MarginRone{1.}
Additionally to its physical significance (and of course related to it), it has given rise to a substantial body of mathematical literature -- possibly even more than the 
\revision{above-mentioned} 
\MarginRone{2.}
classical models.
In this paper, we aim to provide a systematic (though inevitably incomplete) overview 
thereby focusing on \revision{wellposedness} analysis of initial value and 
\revision{time-periodic} 
\MarginRone{28.}
problems, memory and fractional attenuation as well as singular limits and -- with one example each -- control and inverse problems.
}
\section{\revision{Introduction}}
\MarginRone{3.}
The Jordan-Moore-Gibson-Thompson JMGT equation of nonlinear acoustics is a third order in time quasilinear partial differential equation PDE modeling finite amplitude propagation of sound in gases or liquids.
On one hand, nonlinear acoustics has a multitude of applications especially in ultrasonics, ranging from  
high-intensity (focused) ultrasound~\cite{abramov,kennedy2003high, wu2001pathological, yoshizawa2009high} to nonlinear ultrasound tomography~\cite{nonlinparam1, duck2002nonlinear,  nonlinparam3, nonlinparam2, ZHANG20011359, nonlinearity_imaging_JMGT,nonlinearity_imaging_Westervelt}.
On the other hand, this PDE gives rise to interesting mathematical questions.
For this reason, since being put forward by Pedro Jordan~\cite{jordan2008nonlinear,jordan2014second}, 
\MarginRone{a}
the JMGT equation has found 
\revision{enormous} 
\MarginRone{4.}
interest by mathematicians and this has led to a large amount of high quality publications on this model.

The aim of this paper is to review some of the mathematical literature about the JMGT equation. 
In doing so, our strategy is to provide the ideas of a few selected results, and a long (though certainly not complete)
collection of references. Here we mainly focus on the mathematical analysis while clearly making omissions on topics like numerics or 
\revision{multi-physics} 
\MarginRone{6.}
coupling, so as to keep the exposition somewhat concise.
\revision{Nonetheless,}
\MarginRone{5.} 
clearly biased by the scientific interest of the author, a little space has been reserved for highlighting the usefulness of the JMGT model in certain inverse problems.

Distinguishing between two types of nonlinearity inherited from classical Kuznetsov and (Lighthill-)Westervelt models, we will consider
\begin{itemize}
\item 
the JMGT-Westervelt equation for the acoustic pressure $u$
\begin{equation}\label{JMGT-Westervelt}
\tau u_{ttt}+u_{tt}-c^2\Delta u -b \Delta u_t + \eta (u^2)_{tt} +f=0,
\end{equation}
\item 
the JMGT-Kuznetsov equation for the acoustic velocity potential $u$
\begin{equation}\label{JMGT-Kuznetsov}
\tau u_{ttt}+u_{tt}-c^2\Delta u -b \Delta u_t +\bigl(\tilde{\eta}u_t^2+|\nabla u|^2\bigr)_t +f=0,
\end{equation}
\item 
as well as their linearization, usually termed Moore-Gibson-Thompson MGT or Stokes-Moore-Gibson-Thompson \SMGT\, (to acknowledge work on the topic by Stokes~\cite{Stokes}) equation
\begin{equation}\label{MGT}
\tau u_{ttt}+u_{tt}-c^2 \Delta u- b \Delta u_t +f= 0
\end{equation}
\item and its extension by general nonlinearities
\begin{equation}\label{JMGT}
\tau u_{ttt}+u_{tt}-c^2 \Delta u- b \Delta u_t = f(u_t, u_{tt}, \nabla u, \nabla u_{t});
\end{equation}
\end{itemize}

The remainder of this paper is organized as follows.
In Section~\ref{sec:model} we reproduce the derivation of \eqref{JMGT-Westervelt}, \eqref{JMGT-Kuznetsov} from fundamental balance and constitutive equations.
Section~\ref{sec:wellposed} dwells on local and 
\revision{global-in-time} 
\MarginRone{18.}  
 \revision{wellposedness}, also touching upon the question of non-global existence, in the sense of blow-up. In doing so, we consider both settings of posing the PDEs on a bounded domain (along with boundary conditions) and considering it on all of $\mathbb{R}^d$. This includes the addition of nonlocal in time terms to model memory of fractional attenuation. With respect to time, we also look at two options: initial and periodicity conditions.
The limit as $\tau\to0$ is discussed in Section~\ref{sec:tau2zero} and provides a mathematically well-based relation of JMGT to the classical second order in time models of nonlinear acoustics.    
Finally, in Section~\ref{sec:controlinverse} we exemplarily present a control and an inverse problems in the context of the JMGT equation.

While we consider the linear case \eqref{MGT} as a preparatory step for presenting the analysis of the nonlinear one (mainly \eqref{JMGT-Westervelt}), we do not attempt to be exhaustive on the literature about \eqref{MGT} but mainly focus on the nonlinear JMGT setting.

As a minor piece of novelty, we present energy estimates that cover initial and time periodicity conditions for \eqref{JMGT-Westervelt} in a concise and unified manner, thereby also allowing to consider singular limits as $\tau\to0$. As opposed to most of the previous literature, we also treat the nonlinearity completely as a right hand side term, rather than incorporating part of it into the second time derivative term, which makes exposition somewhat easier to follow. 

\section{The model}\label{sec:model}
In a series of papers, Pedro Jordan and co-authors point to the fact that using the Fourier temperature flux law 
\begin{equation} \label{fourier_law}
	\boldsymbol{q}= -\kappa \nabla \theta,
\end{equation}
in the derivation of second-order models of nonlinear acoustics may lead to the so-called paradox of infinite speed of propagation; see~\cite{kuznetsov1971equations, kaltenbacher2009global, kaltenbacher2007numerical, jordan2008nonlinear, jordan2016survey} and 
\revision{suggest using}
\MarginRone{7.} 
the Maxwell--Cattaneo law 
\begin{equation} \label{MC_law}
\boldsymbol{q}+\tau \boldsymbol{q}_t = - \kappa \nabla \theta,
\end{equation}
In here, $\tau>0$ is a time lag, modelling thermal 
\revision{relaxation.}
\MarginRone{8.}
\subsection{Derivation}
\def\ttau{\tilde{\tau}}
We follow the steps taken in~\cite[\S 4.1]{jordan2014second}, see also \cite{fracJMGT}, for a spatially one-dimensional setting. 
The underlying physical assumptions are that the sound wave propagates through a thermally conductive and relaxing liquid or gas with negligible viscosity. 
\MarginRtwo{1.}
\MarginRtwo{2.}
\revision{
Since the acoustic particle velocity $\vec{v}$ is typically irrotational in this context, we can write it as the gradient of a velocity potential $\psi$, i.e., $\vec{v}=\nabla\psi$; in one space dimension, this reads $v=\psi_x$. 
We use the Mach number $\epsilon=\frac{V}{c_0}$ as a small parameter
and introduce the nondimensionalized (nd)
quantities
\[
\begin{aligned}
&\tilde{v}=\frac{v}{V}\text{\ldots nd velocity}
&&\tilde{x}=\frac{x}{L} \text{\ldots nd space}\\
&\tilde{\psi}=\frac{\psi}{LV}\text{\ldots nd velocity potential}
&&\tilde{t}=\frac{t\, c_0}{L} \text{\ldots nd time}\\
&\tilde{p}=\frac{p}{p_0} \text{\ldots nd pressure}
&&
s=\frac{\rho-\rho_0}{\rho_0} \text{\ldots 
condensation}\\
&\tilde{q}=\frac{q}{q_0} \text{\ldots nd heat flux}
&&\tilde{\theta}=\frac{\theta-\theta_0}{\theta_0} \text{\ldots nd temperature}.
\end{aligned}
\]
Here $L$, $V$ $>0$ are the characteristic length and speed;
$p_0$, $q_0$, $\rho_0$, $\theta_0$, $c_0=\sqrt{\frac{p_0}{\rho_0}}$ are the background pressure, heat flux, mass density, temperature, and sound speed, respectively.
This, after skipping the tildes, allows us to write the relevant balance equations 
as 
\begin{equation}\label{mass-momontum}
\begin{aligned}
&s_t+\varepsilon(\psi_x\, s_x+(1+s)\psi_{xx})=0 
&&\text{ mass conservation}\\
&\varepsilon\gamma(1+s)\partial_x\bigl(\psi_t+\tfrac12\varepsilon\,(\psi_x)^2\bigr)+\tilde{p}_x=0 
&&\text{ momontum conservation}
\end{aligned}
\end{equation}
the equation of state, (also called pressure -- density relation) as
\[
p=(1+s)^\gamma\exp(\gamma\mathfrak{e})
\]
with the adiabatic index $\gamma$, and the entropy production law as
\begin{equation}\label{entropyproduction}
\tilde{\kappa}{\mathfrak{e}_t}=-\tilde{K} q_x,
\end{equation}
with $\tilde{\kappa}$ being the dimensionless thermal conductivity, 
$\tilde{K}$ the dimensionless thermal diffusivity, and $\mathfrak{e}$ the dimensionless entropy.
We will also use $\ttau$ as a dimensionless version of the relaxation time, the so-called Knudsen number.
Employing a weakly-nonlinear approximation
\begin{equation} \label{weakly_nl_assumptions_alpha}
	\epsilon <<1,\quad s=O(\epsilon),\quad \theta=O(\epsilon),\quad \tilde{K}=O(\epsilon), \quad \ttau=O(\epsilon), \quad |\mathfrak{e}|=O(\epsilon^2)
\end{equation}
allows us to expand $p$ as
\[
p=1+\gamma\, s + \tfrac12\gamma(\gamma-1)\, s^2 + \gamma \mathfrak{e}+O(\varepsilon^3).
\]
Differentiating this expression with respect to space and inserting it into the momentum balance equation (second line of \eqref{mass-momontum}), after multiplication with $\frac{1}{1+s}=1-s+O(\varepsilon^2)$ and skipping $O(\varepsilon^2)$ terms yields
\[
\partial_x\Bigl(\psi_t+\tfrac12\epsilon(\psi_x)^2+\varepsilon^{-1}\bigl(s+\tfrac12(\gamma-2) s^2+\mathfrak{e}\bigr)\Bigr)=0.
\] 
Integrating with respect to space and differentiating with respect to time we obtain
\[
\psi_{tt}+\tfrac12\epsilon\partial_t(\psi_x)^2
+\varepsilon^{-1}\bigl(1+(\gamma-2)s\bigr)\, s_t
+\varepsilon^{-1}\mathfrak{e}_t=0.
\] 
Now we use mass balance (first line of \eqref{mass-momontum}) to eliminate $s_t$, which 
}
leads to the equation
\begin{equation} \label{eq_0}
	\begin{aligned}
		\psi_{tt}+\tfrac12 \epsilon \partial_t(\psi_x)^2-(1+(\gamma-2)s)[\psi_x s_x+(1+s)\psi_{xx}]=-\epsilon^{-1}\mathfrak{e}_t,
	\end{aligned}
\end{equation}
for the acoustic velocity potential $\psi$, 
see~\cite[(53)]{jordan2014second}. 
The alternative heat flux law \eqref{MC_law} in a dimensionless spatially 1-d version
\[
(1+\ttau \partial_t) q(t) = -\tilde{\kappa} \theta_x,
\]
comes into play via the entropy production law \eqref{entropyproduction},
whose combination with the above line 
yields the following entropy equation
\begin{equation}\label{entropyeq_0}   
(1+\ttau  \partial_t) \mathfrak{e}_t=\tilde{K}  \theta_{xx}.
\end{equation}
Utilizing the approximations $s = -\epsilon \psi_t+ O(\epsilon^2)$ cf. \cite[(49)]{jordan2014second} and $\theta = - \epsilon (\gamma-1) \psi_t+ O(\epsilon^2)$ \cite[(57),(58)]{jordan2014second}, neglecting all $O(\epsilon^2)$ terms, we rewrite \eqref{eq_0} and \eqref{entropyeq_0} as  
\begin{equation} \label{eq_1}
	\begin{aligned}
		\psi_{tt}+\tfrac12 \epsilon \partial_t(\psi_x)^2-(1-(\gamma-2)\epsilon\psi_t)[-\epsilon \psi_x \psi_{tx}+(1-\epsilon \psi_t)\psi_{xx}]=-\epsilon^{-1}\mathfrak{e}_t;
	\end{aligned}
\end{equation}
and
\begin{equation} \label{e_t_GFE}
	\begin{aligned}
		(1+\ttau  \partial_t)\mathfrak{e}_t=-\epsilon\tilde{K}(\gamma-1)\psi_{txx};
	\end{aligned}
\end{equation}
Elimination of $\mathfrak{e}$ can be achieved by applying $(1+\ttau \partial_t)$ to \eqref{eq_1} and dividing \eqref{e_t_GFE} by $\epsilon$, which leads to
\begin{equation} \label{eq_2_I}
	\begin{aligned}
(1+\ttau  \partial_t)\left\{\psi_{tt}+\tfrac12 \epsilon \partial_t(\psi_x)^2-(1-(\gamma-2)\epsilon \psi_t)[-\epsilon \psi_x \psi_{tx}+(1-\epsilon \psi_t)\psi_{xx}]\right\}\\
			= \,  \tilde{K}(\gamma-1) \psi_{txx}. 
	\end{aligned}
\end{equation}
Since $\ttau =O(\epsilon)$, neglecting the $O(\epsilon^2)$ terms in \eqref{eq_2_I} yields
\begin{equation} \label{eq_3}
	\begin{aligned}
		(1+\ttau  \partial_t)\psi_{tt}-\ttau \partial_t \psi_{xx}-(1-\epsilon (\gamma-1)\psi_t)\psi_{xx}+\epsilon \partial_t(\psi_x)^2= \tilde{K}(\gamma-1) \psi_{txx}.
	\end{aligned}
\end{equation}
Finally, we multiply with $(1-\epsilon (\gamma-1)\psi_t)^{-1} = 1+\epsilon (\gamma-1)\psi_t + O(\epsilon^2)$ and neglect all $O(\epsilon^2)$ terms to obtain
\begin{equation} \label{final_eq_2}
\ttau  \psi_{ttt}+(1+\epsilon (\gamma-1)\psi_t)\psi_{tt}-\psi_{xx}-\ttau \psi_{txx}
-\tilde{K}(\gamma-1) \psi_{txx}+\epsilon \partial_t(\psi_x)^2= 0. 
\end{equation}
This obviously is a dimensionless 1-d version of \eqref{JMGT-Kuznetsov}.
\revision{Neglecting} 
\MarginRone{9.}
non-cumulative nonlinear effects by 
\revision{employing} 
\MarginRone{10.}
the approximation $|\nabla \psi|^2 \approx c^{-2}\psi_t^2$, leads to \eqref{JMGT-Westervelt}.

%
\section{Well-posedness analysis}\label{sec:wellposed}
In what follows, we aim to provide a unified exposition over the various topics, relying on a certain energy identity, from which local and, in the dissipative case $b>c^2\tau$ also global \revision{wellposedness} as well as exponential decay, but also existence of 
\revision{time-periodic} 
\MarginRone{28.}
solutions and limits as $\tau\to0$ can be derived for \eqref{JMGT-Westervelt}. The analysis of \eqref{JMGT-Kuznetsov} follows similar lines but require sophisticated higher order energy estimates, which we do not 
\revision{explicitly} 
\MarginRone{11.}
show here.

\subsection{Linear(ized) versions}\label{sec:lin}
We start by recalling the 
\revision{\SMGT\, equation}
\MarginRone{12.}
\begin{equation}\label{MGT_ana}
	\begin{aligned}
		\tau u_{ttt}+u_{tt}-c^2 \Delta u - b \Delta u_t = 0.
	\end{aligned}
\end{equation}
\revision{As already pointed out in \cite{jordan2014second}, 
it
arises as a model in several contexts outside acoustics; in particular \cite{dell2017moore,pellicer2019optimal} detail the relation to a standard model of linear viscoelesticity.}
\MarginRone{b}

The key sufficient criterion 
\begin{equation}\label{delta}
\delta:=b-\tau c^2>0
\end{equation}
for (exponential) stability of this linear system can be nicely motivated by the Routh-Hurwitz stability criterion from control theory. 
To this end, we follow the exposition in \cite[Section 1.1]{RackeSaidHouari:2021}.
Considering \eqref{MGT_ana} either 
\begin{itemize}
\item[(a)] on a bounded domain $\Omega\subset\mathbb{R}^d$ and equipping $-\Delta$ with boundary conditions that render its inverse $(-\Delta_B)^{-1}$ a positive definite 
\revision{self-adjoint} 
\MarginRone{13.}
compact operator on $L^2(\Omega)$ or 
\item[(b)] on all of $\mathbb{R}^d$, 
\end{itemize}
we can use (a) an eigenfunction expansion of the negative Laplacian or (b) apply the Fourier transform $\mathcal{F}_x$ with respect to space to arrive at the family of ODEs
\[
\tau \tilde{u}'''+\tilde{u}''+c^2 \zeta \tilde{u} + b \zeta \tilde{u}' = \tilde{f}  
\]   
with (a) $\zeta=\lambda_j$, $\tilde{u}(t)=\langle u(t),\varphi^j\rangle$, $\lambda_j$, $\varphi^j$ eigenvalues and -functions of $-\Delta_B$ or (b) $\zeta=|\xi|^2$, $\tilde{u}(t)=(\mathcal{F}_x u(t))(\xi)$, $\xi\in \mathbb{R}^d$.
This can be written as a first order in time system 
\[
U'=Au+F\quad \text{ with }A=\left(\begin{array}{ccc}
0&1&0\\
0&0&1\\
-1&-b\zeta/\tau &-c^2\zeta/\tau
\end{array}\right)
\]
The Routh-Hurwitz stability criterion now states that the real parts of the eigenvalues of $A$ have negative real part (implying stability of the system) iff the principal minors of the Hurwitz matrix are positive
\[
\sigma_p(A) \subseteq\mathbb{R}^-+\imath\mathbb{R}
\quad \Leftrightarrow \quad 
\left(m_j
\left(\begin{array}{ccc}
1&\tau&0\\
c^2\zeta&b\zeta&1\\
0&0 &c^2\zeta
\end{array}\right)>0, \ j\in\{1,2,3\}\right).
\]
Since these minors compute as $m_1=1$, $m_2=(b-\tau c^2)\zeta$, $m_3=c^2\zeta\,m_2$, equivalence to \eqref{delta} is obvious.

This elementary observation is confirmed and refined by semigroup methods \cite{kaltenbacher2011wellposedness,marchand2012abstract}, according to which \eqref{MGT}
generates a continuous semigroup (in fact, a group) which is exponentially stable under condition \eqref{delta}, with a decay factor that has been quantified in \cite{pellicer2019optimal} by proving normality of the generator with respect to an appropriately chosen inner product.
The so-called critical (or inviscid) case $\delta=b-\tau c^2=0$ leads to marginal stability, while $\delta<0$ according to \cite{conejero2015chaotic} leads to chaotic behaviour.
In \cite{BucciEller2021} the authors point to the hyperbolic nature of the equation, as well as the fact that the spatial boundary is characteristic. 

\medskip

To tackle the quasilinear case by means of 
\revision{fixed-point} 
\MarginRone{14.}
theorems, we need to consider 
linearization of \eqref{JMGT-Westervelt}, \eqref{JMGT-Kuznetsov}, 
\begin{equation}\label{JMGT_lin}
	\begin{aligned}
		&\tau u_{ttt}+u_{tt}-c^2 \Delta u- b \Delta u_t = f 
	\end{aligned}
\end{equation}
\MarginRone{c}
with $f=-\eta (u^2)_{tt}$ for \eqref{JMGT-Westervelt} 
and $f=-\bigl(\tilde{\eta}u_t^2+|\nabla u|^2\bigr)_t$ for \eqref{JMGT-Kuznetsov}.

While we have the physical case $\Omega\subseteq\mathbb{R}^d$, $d=3$ in mind throughout this paper, we remain with a simpler setting with respect to boundary conditions (Dirichlet) and nonlinearity (\eqref{JMGT-Westervelt} rather than \eqref{JMGT-Kuznetsov}) to keep exposition transparent here, while pointing to the more general situations that can be found in the cited literature.


\subsection{Energy estimates}\label{subsec:energyestimates}
As a preparation for handling the nonlinearity, we here showcase energy estimates, as they also give some intuition with a minimal amount of mathematical machinery. 
To this end, we assume \eqref{JMGT_lin} to hold on either (a) a bounded domain $\Omega\subseteq\mathbb{R}^d$ with homogeneous Dirichlet boundary conditions or (b) all of $\mathbb{R}^d$.

The testing strategy that predominates in this scenario is motivated by the fact that the combined quantity 
\begin{equation}\label{defz}
z:=\tau u_t+u 
\end{equation}
solves a perturbed 
wave equation
\begin{equation}\label{wave_z}
z_{tt}-(c^2+\tfrac{\delta}{\tau})\Delta z 
= f
-\tfrac{\delta}{\tau}\Delta u 
\end{equation}
\MarginRone{c}
and the standard test function for deriving energy estimate for the wave equation is the first time derivative of the state $z_t$. 
\revision{However, the resulting regularity $z\in W^{1,\infty}(0;T;L^2(\Omega))\cap L^\infty(0,T;H^1(\Omega))$ would not suffice for controlling the norm of the nonlinear term that is needed to close the estimates.
To obtain sufficient spatial regularity, we therefore 
}
\MarginRtwo{3.}
\revision{test with 
$- \Delta z_t=-\Delta(\tau u_{tt}+u_t)$, 
augmented by a small term which}
\MarginRone{15.}  
provides an energy contribution that is zero order in time in $u$. 
That is, we multiply \eqref{JMGT_lin} with $-\Delta(\tau u_{tt}+\sigma u_t+\rho u)$ with $\sigma>0$ close to one and $\rho>0$ small enough. 
Integrating (by parts) over space and time and using the identities
\[
\begin{aligned}   
&\prodLtwo{\Delta u}{\Delta u_{tt}}=\frac{d}{dt}\prodLtwo{\Delta u}{\Delta u_t} -\nLtwo{\Delta u_t}^2\\
&\prodLtwo{\nabla u_{ttt}}{\nabla u_t}=\frac{d}{dt}\prodLtwo{\nabla u_{tt}}{\nabla u_t}-\nLtwo{\nabla u_{tt}}^2\\
&\prodLtwo{\nabla u_{ttt}}{\nabla u}=\frac{d}{dt}\prodLtwo{\nabla u_{tt}}{\nabla u}-\frac{d}{dt}\frac12\nLtwo{\nabla u_t}^2\\
&\prodLtwo{\nabla u_{tt}}{\nabla u}=\frac{d}{dt}\prodLtwo{\nabla u_t}{\nabla u}-\nLtwo{\nabla u_t}^2
\end{aligned}
\]
we obtain the energy identity
\begin{equation}\label{enid}
\begin{aligned}   
&\frac12 \tau^2\nLtwo{\nabla u_{tt}}^2\Big|_0^t+\tau(1-\sigma)\int_0^t\nLtwo{\nabla u_{tt}}^2\, ds\\
&+\frac12 \tau b\nLtwo{\Delta u_t}^2\Big|_0^t+(b\sigma-\tau c^2)\int_0^t\nLtwo{\Delta u_t}^2\, ds\\
&+\frac12 (\sigma c^2+b\rho)\nLtwo{\Delta u}^2\Big|_0^t+c^2\rho\int_0^t\nLtwo{\Delta u}^2\, ds\\
&+\frac12 (\sigma-\tau\rho)\nLtwo{\nabla u_t}^2\Big|_0^t-\rho\int_0^t\nLtwo{\nabla u_t}^2\, ds\\
&+\tau c^2\prodLtwo{\Delta u_t}{\Delta u}\Big|_0^t+\tau\sigma\prodLtwo{\nabla u_{tt}}{\nabla u_t}\Big|_0^t\\
&+\tau\rho\prodLtwo{\nabla u_{tt}}{\nabla u}\Big|_0^t + \rho\prodLtwo{\nabla u_t}{\nabla u}\Big|_0^t\\
&=\int_0^t\prodLtwo{f}{-\Delta(\tau u_{tt}+\sigma u_t+\rho u)}.
\end{aligned}
\end{equation} 
Considering the 
\revision{leading-order} 
\MarginRone{16.}  
contributions, that is, the first terms in the first and second line of \eqref{enid}, it gets apparent that the conditions 
\begin{equation}\label{condcoeff_local}
\tau>0, b=\tau c^2 +\delta>0 
\end{equation}
are sufficient for enabling an energy estimate that allows to bound a solution $u$, provided $f$ is regular enough, since all other terms containing $u$ can be dominated by either 
$\nLtwo{\nabla u_{tt}(t)}^2$ or $\nLtwo{\Delta u_{t}(t)}^2$ in a Gronwall type argument. 
This is the basis for a proof of 
\revision{local-in-time} 
\MarginRone{17.}  
\revision{wellposedness} of the nonlinear problem.
In order to show 
\revision{global-in-time} 
\MarginRone{18.}  
 \revision{wellposedness} and exponential decay, we require positivity of all terms containing norms of derivatives of $u$, that is the first seven terms on the \revision{left-hand-side} \MarginRone{19.}, as the eighth one can be controlled by the estimate
\begin{equation}\label{PF}
\nLtwo{\nabla v}\leq \|\nabla v\|_{H^1(\Omega)}\leq C(\Omega)\nLtwo{\Delta v} \quad v\in H^2(\Omega)\cap H_0^1(\Omega)
\end{equation}
that follows from  elliptic regularity.
To this end, we additionally assume  
\begin{equation}\label{condcoeff_global}
\delta=b-\tau c^2>0 
\end{equation}
and choose $\sigma:=1-\min\{\frac{\delta}{2b}, \frac{\tau\delta}{c^2}\}<1$, which makes the first seven terms on the \revision{left-hand-side} \MarginRone{19.} of \eqref{enid} positive and additionally allows to control the 
\revision{ninth}
\MarginRone{20.} 
and tenth term by some of the previous norms 
\[
\begin{aligned}
&\tau c^2\prodLtwo{\Delta u_t}{\Delta u}
<\frac12 \tau b\nLtwo{\Delta u_t}^2+\frac12 (\sigma c^2+b\rho)\nLtwo{\Delta u}^2\\
&\tau\sigma\prodLtwo{\nabla u_{tt}}{\nabla u_t}
< \frac12 \tau^2\nLtwo{\nabla u_{tt}}^2+\frac12 (\sigma-\tau\rho)\nLtwo{\nabla u_t}^2
\end{aligned}
\] 
for $\rho>0$ small enough. All remaining terms are lower order and can be dominated by higher order ones due to \eqref{PF}, by choosing $\rho>0$ small enough.
This leads us to defining an energy by 
\begin{equation} \label{def_energy}
\begin{aligned}
\mathcal{E}[u](t):=&
\tau^2\nLtwo{\nabla u_{tt}(t)}^2
+\tau\nLtwo{\Delta u_t(t)}^2\\
&+\nLtwo{\nabla u_t(t)}^2
+\nLtwo{\Delta u(t)}^2.
\end{aligned}
\end{equation}

The right hand side can then be estimated by means of Young's inequality, weighting the terms containing $u$ in such a way that they can be dominated by \revision{left-hand-side} \MarginRone{19.} terms
\begin{equation}\label{boundf}
\begin{aligned}
&\int_0^t\prodLtwo{f}{-\Delta(\tau u_{tt}+\sigma u_t+\rho u)}\, ds
=\int_0^t\prodLtwo{\nabla f}{\nabla(\tau u_{tt}+\sigma u_t+\rho u)}\, ds\\
&\leq \frac{\epsilon}{2}\int_0^t
\mathcal{E}[u](s)\, ds
+ \frac{1}{2\epsilon}\Bigl(1+\sigma^2+\rho^2\Bigr)\int_0^t\nLtwo{\nabla f}^2\, ds,
\end{aligned}
\end{equation}
where we have assumed that $f$ satisfies homogeneous 
\revision{Dirichlet}  
\MarginRone{21.} 
boundary condition to allow for integration by parts without adding boundary terms.

\revision{
Altogether, we obtain the following energy estimate.
\begin{proposition}
For the energy defined by \eqref{def_energy} 
\begin{itemize}
\item 
under the condition $b=\tau c^2+\delta>0$, we have
\begin{equation} \label{energy_est}
\mathcal{E}[u](t) \leq C_0\Bigl(\mathcal{E}[u](0)+\int_0^t\bigl(\mathcal{E}[u](s) +\nLtwo{\nabla f}^2\bigr)\, ds \Bigr);
\end{equation}
\item 
in case $\delta=b-\tau c^2>0$, we even have  
\begin{equation} \label{energy_est_dissip}
\mathcal{E}[u](t) +c_1\int_0^t\mathcal{E}[u](s)\,ds\leq C_0\Bigl(\mathcal{E}[u](0)+\int_0^t\nLtwo{\nabla f}^2\, ds\Bigr).
\end{equation}
\end{itemize}
\end{proposition}
The dissipation present in \eqref{energy_est_dissip} under condition $\delta=b-\tau c^2>0$}
\MarginRtwo{3.}
helps us to establish 
\revision{global-in-time} 
\MarginRone{22.}  
wellposedness of the nonlinear problem.
Note that the small factor $\epsilon>0$ in \eqref{boundf} can be chosen independently of $\tau$, but for establishing \eqref{energy_est_dissip} it must be adapted to $\delta>0$ since the dissipative terms in \eqref{enid} depend on $\delta$. Thus, when considering parameter limits, one has to take into account the fact that the constants $C_0$, $c_1$ in \eqref{energy_est_dissip} are independent of $\tau$ but may depend on $\delta$. 
This conforms to the fact that limits as $\delta\searrow0$ can only be established in the context of 
\revision{local-in-time} 
\MarginRone{17.}  
\revision{wellposedness} results for \eqref{JMGT-Westervelt}, cf. \cite{b2zeroJMGT}.
\subsection{Initial value problem with small data}\label{sec:ivp}
In order to present the ideas, we focus on the JMGT-Westervelt equation \eqref{JMGT-Westervelt} on a smooth bounded domain $\Omega\subseteq\mathbb{R}^d$ with homogeneous Dirichlet conditions to maximize comparability with results in the literature and minimize the amount of technicalities (as \eqref{JMGT-Kuznetsov} requires higher order energy estimates). 
We consider the initial boundary value problem
\begin{equation}\label{JMGT_ivp}
	\begin{aligned}
		&\tau u_{ttt}+u_{tt}-c^2 \Delta u- b \Delta u_t = -\eta(u^2)_{tt} &&\text{ in }(0,T)\times\Omega\\
		&u=0 &&\text{ on }(0,T)\times\partial\Omega\\
	    &u(0)=u_0, \quad u_t(0)=u_1, \quad u_{tt}(0)=u_2		
	\end{aligned}
\end{equation}
and its linear version
\begin{equation}\label{JMGT_lin_ivp}
	\begin{aligned}
		&\tau u_{ttt}+u_{tt}-c^2 \Delta u- b \Delta u_t = f &&\text{ in }(0,T)\times\Omega\\
		&u=0 &&\text{ on }(0,T)\times\partial\Omega\\
	    &u(0)=u_0, \quad u_t(0)=u_1, \quad u_{tt}(0)=u_2.	
	\end{aligned}
\end{equation}
The initial conditions are supposed to be chosen such that 
\begin{equation} \label{initial_energy}
\mathcal{E}[u](0)=\tau^2\nLtwo{\nabla u_2}^2
+\tau\nLtwo{\Delta u_1}^2+\nLtwo{\nabla u_1}^2
+\nLtwo{\Delta u_0}^2<\infty.
\end{equation}
A smallness assumption on  $\mathcal{E}[u](0)$ is needed to prove 
\revision{global} 
\MarginRone{d}
existence of solutions to the initial value problem. Indeed, blow-up in finite time can be shown otherwise, cf. \cite{NikolicWinkler2024_blow-up} which is discussed in Section~\ref{sec:blowup} below.
\subsubsection{Local-in-time \revision{wellposedness}}\label{subsec:local}
In this subsection, we fix the time horizon $T>0$ and assume that $b=c^2\tau+\delta>0$ so that \eqref{energy_est} holds.
The typical 
\revision{local-in-time} 
\MarginRone{17.}  
 \revision{wellposedness} proofs found in the literature on the JMGT equation rely on Banach's Contraction Principle for the 
\revision{fixed-point} 
\MarginRone{14.}
operator 
\MarginRtwo{4.}
\[
\mathcal{T}:\revision{v}\mapsto u \text{ solving \eqref{JMGT_lin_ivp} with }f=-\eta(\revision{v}^2)_{tt}=-2\eta(\revision{v}\revision{v_{tt}}+\revision{v_t}^2). 
\] 
We will here focus on the verification of a self-mapping property on a sufficiently small ball (with respect to the norm induced by the energy \eqref{def_energy} as this is probably the most transparent and illustrative part of the proof.  
The key ingredients for this purpose are the energy estimate \eqref{energy_est} and a bound on $\int_0^t\nLtwo{\nabla f}^2\, ds$ in terms of $\mathcal{E}[\revision{v}]$. Note that $f$ inherits homogeneous Dirichlet boundary conditions from $\revision{v}$ and thus satisfies the conditions underlying \eqref{boundf}.
\begin{equation}\label{estf}
\begin{aligned}
&\nLtwo{\nabla f}^2\, ds = 2\eta \nLtwo{\revision{v_{tt}}\nabla \revision{v} + \revision{v}\nabla \revision{v_{tt}}+2\revision{v_t}\nabla \revision{v_t}}^2\\
&\leq 8\eta\Bigl(\|\revision{v_{tt}}\|_{L^4(\Omega)}^2\|\nabla \revision{v}\|_{L^4(\Omega)}^2 
+ \|\revision{v}\|_{L^\infty(\Omega)}^2 \nLtwo{\nabla \revision{v_{tt}}}^2
\\&\qquad\qquad
+ 2\|\revision{v_t}\|_{L^\infty(\Omega)}^2 \nLtwo{\nabla \revision{v_t}}^2\Bigr)\\
&\leq C (1+\tau^{-2}) \mathcal{E}[\revision{v}]^2,
\end{aligned}
\end{equation}
due to Sobolev embeddings and elliptic regularity.
Here we assume $\eta$ to be constant, but corresponding results with space dependent coefficients are possible as well \cite{nonlinearity_imaging_JMGT}.
The constant $C$ does not depend on the time horizon nor on $\tau$. A refined analysis allows to eliminate dependence on $\tau$, thus enabling a limiting analysis as $\tau\searrow0$ \cite{bongarti2020vanishing,JMGT_BKVN,JMGT_Neumann}.

Gronwall's inequality allows us to conclude from \eqref{energy_est}, \eqref{estf} that $\mathcal{T}$ is a self-mapping on the set 
\[
\begin{aligned}
\mathcal{W}=\{w\in L^2(0,T;H_0^1(\Omega))\, : \,& \sup_{t\in(0,T)} \mathcal{E}[w](t)\leq \rho, \\ 
&(w(0),w_t(0),w_{tt}(0))=(u_0,u_1,u_2)\}
\end{aligned}
\]
as long as $\mathcal{E}[u](0)\leq \rho_0$ and $\rho$, $\rho_0$ are small enough.
\revision{A contractivity proof  relies on the fact that the difference between values $\mathcal{T}(v_1)-\mathcal{T}(v_2)$ satisfies \eqref{JMGT_lin_ivp} with 
$f=-\eta\bigl((v_1+v_2)(v_1-v_2)\bigr)_{tt}
=-\eta\bigl((v_{1,tt}+v_{2,tt})(v_1-v_2)+2(v_{1,t}+v_{2,t})(v_{1,t}-v_{2,t})+(v_1+v_2)(v_{1,tt}-v_{2,tt})\bigr)$, 
which can be estimated analogously to \eqref{estf}.
\footnote{Note that we can here rely on the same topology for self-mapping and contractivity, while several contractivity or uniqueness proofs in other models of nonlinear acoustics resort to a weaker topology.}
}
\MarginRtwo{5.}
Combining this with 
\revision{the self-mapping property shown above,}
\revision{local-in-time} 
\MarginRone{17.}  
\revision{wellposedness} can be shown for \eqref{JMGT-Westervelt} on a smooth bounded domain $\Omega$ in case $b=c^2\tau+\delta>0$ with small enough initial data \cite{kaltenbacher2012well,JMGT_BKVN,NikolicWinkler2024_blow-up}.
%
\subsubsection{Global-in-time \revision{wellposedness} and exponential decay}\label{subsec:global}
In the dissipative case $\delta=b-\tau c^2>0$, a combination of \eqref{energy_est_dissip}  with \eqref{estf} allows to deduce the estimate 
\begin{equation} \label{energy_est_barrier}
\mathcal{E}[u](t) +c_1\int_0^t\mathcal{E}[u](s)\,ds\leq C_1\Bigl(\mathcal{E}[u](0)+\int_0^t\mathcal{E}[u](s)^2\, ds\Bigr)
\end{equation}
for a solution $u$ to \eqref{JMGT_ivp} on the time interval $(0,T)$, where both constants $c_1$, \revisionown{$C_1\geq1$} are independent of $T$. This is the basis for a 
\revision{global-in-time} 
\MarginRone{22.}  
\revision{well-posedness}
\MarginRone{23.}  
proof for sufficiently small initial data 
as follows.
\revision{Estimate \eqref{energy_est_barrier} implies
\[
\begin{aligned}
&m(t)\leq C_1\mathcal{E}[u](0)+ \frac{C_1}{c_1} m(t)^2\\
&\text{ for }m(t):=\sup_{s\in(0,t)}\mathcal{E}[u](s) +c_1\int_0^t\mathcal{E}[u](s)\,ds.
\end{aligned}
\]
Therefore, we can apply the following result.
\begin{quote}
\cite[Lemma 3.7]{Strauss1968}
Let $m(t)$ be a non-negative continuous function of $t$ satisfying the inequality 
\[
m(t)\leq \hat{c}_1+\hat{c}_2 m(t)^\gamma
\]
in some interval containing $0$, where $\hat{c}_1$ and $\hat{c}_2$ are positive constants and $\gamma>1$. 
If $m(0)\leq \hat{c}_1$ 
and 
$\hat{c}_1\, \hat{c}_2^{(\gamma-1)^{-1}} < (1-\gamma^{-1})\,\gamma^{-(\gamma-1)^{-1}}$, 
then in the same interval
\[
m(t)<\frac{\hat{c}_1}{1-\gamma^{-1}}.
\]
\end{quote}
Indeed, setting $\hat{c}_1:=C_1\mathcal{E}[u](0)$, $\hat{c}_2:=\frac{C_1}{c_1}$, $\gamma=2$, we deduce
\[
\mathcal{E}[u](t) +c_1\int_0^t\mathcal{E}[u](s)\,ds\leq 2 C_1\mathcal{E}[u](0),
\]
provided $\mathcal{E}[u](0)<\frac{c_1}{C_1^2}$.
This implies that on any interval $(0,T)$, the energy $\mathcal{E}[u](t)$ can be bounded by a constant independent of $T$ and therefore $u$ exists (and its energy remains bounded) on all of $(0,\infty)$. 
}
\MarginRone{e}
In fact, it is readily checked that also the differentiated version
\[
0\geq \frac{d}{dt}\mathcal{E}[u](t) + \revisionown{c_1}\mathcal{E}[u](t)
=e^{-\revisionown{c_1}t}\, \frac{d}{dt}\Bigl(e^{\revisionown{c_1} t}\mathcal{E}[u](t)\Bigr)
\quad t>0
\]
holds true. From this, we can immediately conclude exponential decay as it implies
\[
e^{\revisionown{c_1} t}\mathcal{E}[u](t)\leq e^{\revisionown{c_1} 0}\mathcal{E}[u](0)=\mathcal{E}[u](0)
\quad t>0. 
\]

In \cite{BongartiLasiecka2022} it is shown that global \revision{wellposedness} \MarginRone{23.}and exponential decay of solutions to \eqref{JMGT-Westervelt} can even be achieved in the critical case $\delta=0$ by an appropriate boundary feedback -- so-called absorbing boundary conditions, see also \cite{bucci2018feedback} and \eqref{MGT-Neumann} in Section~\ref{subsec:control} below.
\subsubsection{Cauchy problem}
Considering \eqref{JMGT-Westervelt} or \eqref{JMGT-Kuznetsov} 
\revision{
on all of $\mathbb{R}^d$, rather than on a bounded domain, leads to significant differences,
particularly with respect to the provable decay rates.} 
\MarginRone{25.}
Like in the previous subsection, we consider the dissipative setting $\delta=b-\tau c^2>0$ here.
First of all, for the linear \SMGT\, equation, in \cite{PellicerSaidHouari:2019} the pointwise ODE resulting from application of the spatial Fourier transform to \eqref{MGT}
\begin{equation}\label{MGT_Fourier}
\tau \hat{u}_{ttt}+\hat{u}_{tt}+c^2|\xi|^2\hat{u}+ b |\xi|^2 \hat{u}_t = 0
\end{equation}
is considered.
Using an appropriately constructed Lyapunov function, a pointwise estimate in spatial Fourier domain of the form
\[
\begin{aligned}
&|\hat{V}(\xi,t)|^2:=|\tau\hat{u}_{tt}(\xi,t)+\hat{u}_t(\xi,t)|^2+|\tau\nabla\hat{u}_t(\xi,t)+\nabla\hat{u}(\xi,t)|^2+|\nabla\hat{u}_t(\xi,t)|^2\\
&\leq C\exp\Bigl(-c\frac{|\xi|^2}{|\xi|^2+1}t\Bigr)|\hat{V}(\xi,0)|^2
\end{aligned}
\]
is established by Pellicer and Said-Houari in \cite{PellicerSaidHouari:2019}, which allows them to prove the decay estimate
\[
\|\nabla^j V(t)\|_{L^2(\mathbb{R}^d)}^2\leq C (1+t)^{-\frac{d}{4}-\frac{j}{2}} \|V(0)\|_{L^1(\mathbb{R}^d)}^2
+ C e^{-ct} \|\nabla^j V(0)\|_{L^2(\mathbb{R}^d)}^2,
\]
provided the right hand side is finite. 
It is in fact the low frequency part $\{\xi\in \mathbb{R}^d\, : \, |\xi|\leq 1\}$ that slows down decay as compared to the bounded domain setting.
Note that the definition of $V$ contains the physical wave energy of the combined quantity $z$ defined by \eqref{defz} and satisfying the second order wave equation \eqref{wave_z}.
Global \revision{wellposedness} of the Cauchy problem and a corresponding decay result is proven by Racke and Said-Houari \cite{RackeSaidHouari:2021} in the fully nonlinear setting of \eqref{JMGT-Kuznetsov}, thus including even the gradient nonlinearity, which requires sophisticated estimates in higher order norms, see \cite{RackeSaidHouari:2021}.
\subsection{Blow up in finite time}\label{sec:blowup}
\def\Tmax{T_{\text{max}}}
\def\limsup{\text{lim sup}}
An important question complementary to global-in-time \revision{wellposedness} for small initial data is whether nonexistence of global solutions with large initial data can be proven and how solutions behave near the end of the (finite) maximal existence time horizon $T_{\text{max}}$. 
In \cite{NikolicWinkler2024_blow-up}, Nikoli\'{c} and Winkler study a general quasilinear JMGT-type model comprising \eqref{JMGT-Westervelt} 
\begin{equation}\label{JMGT_gutt}
\tau u_{ttt}+\alpha u_{tt}-c^2 \Delta u- b \Delta u_t = g(u)_{tt}
\end{equation}
that is \eqref{JMGT} with 
$f(u_t, u_{tt}, \nabla u, \nabla u_{t})
=g'(u)u_{tt}+g''(u)u_t^2$
under the assumptions $\tau>0$, $c^2>0$, $b>0$, $\alpha\in\mathbb{R}$, 
$g(0)=0$, $g\in C^3(\mathbb{R})$, comprising \eqref{JMGT-Westervelt}.
They prove that for a strong solution of \eqref{JMGT_gutt} on a smooth bounded domain $\Omega\subseteq\mathbb{R}^3$ with homogeneous Dirichlet boundary conditions and arbitrary regular enough initial data, the implication
\[
\Tmax<\infty \quad \Longrightarrow \quad \limsup_{t\nearrow\Tmax}\|u(t)\|_{L^\infty(\Omega)}=\infty
\]
holds \cite[Theorem 1.1]{NikolicWinkler2024_blow-up}. Under the additional assumptions that the first eigenfunction of the Dirichlet Laplacian on $\Omega$ is strictly positive and $g$ satisfies
\[
g''\geq0\, \quad 
\lim_{\xi\to\infty}\frac{g(\xi)}{\xi}=\infty, \quad 
\int_{\xi_0}^\infty\frac{1}{g(\xi)}\, d\xi <\infty \text{ for some }\xi_0>0
\]
it is shown that indeed $\Tmax<\infty$ must hold for any initial data whose projection on the first Dirichlet eigenfunction is large enough \cite[Theorem 1.2]{NikolicWinkler2024_blow-up}.
The proof is carried out by excluding gradient blow-up, that is, proving that if on the contrary the $L^\infty(\Omega)$ norm of $u(t)$ remains bounded as $t$ tends to a finite $\Tmax$, then also the $L^2(\Omega)$ norms of $\Delta u(t)$, $\nabla u_t(t)$, $\Delta u_t(t)$, $\nabla u_{tt}(t)$ must stay bounded.
The additional regularity requirements on the initial data as compared to finiteness of the energy \eqref{initial_energy} (whose smallness is required for the ``opposite'' result of 
\revision{global-in-time} 
\MarginRone{22.}  
 \revision{wellposedness}) is 
$u(0)\in H^4(\Omega)$, $u_t(0)\in H^2(\Omega)$, $u_{tt}(0)\in H^2(\Omega)$.
\\
We also refer to \cite{chen2019nonexistence} for nonexistence results in the semilinear setting \eqref{JMGT} with $f(u_t, u_{tt}, \nabla u, \nabla u_{t})$ replaced by $|u|^p$ in a certain range of powers $p$, to \cite{ChenPalmieri:2021derivative} for blow-up with $f(u_t, u_{tt}, \nabla u, \nabla u_{t})=|u_t|^p$, with $1<p<\frac{d+1}{d-1}$, and to \cite{MingYangFanYao2021} for a combination replacing $f(u_t, u_{tt}, \nabla u, \nabla u_{t})$ by $|u|^p+|u_t|^p$.  
\subsection{JMGT with memory and fractional attenuation}
Probably the largest amount of work in the literature in the context of the (J)MGT equation 
has been dedicated to studying the addition of memory. As we do not feel able to give proper credit to all this work, we confine ourselves to only highlighting 
the pioneering results from \cite{LasieckaWang15a,LasieckaWang15b}
and providing a 
(probably still incomplete) 
list of further references.

To start with, let us recall the fact that the (J)MGT equation itself can be viewed as a wave equation with memory for the combined quantity $z:=\tau u_t+u$ cf. \eqref{defz}, that due to \eqref{wave_z} and resolving \eqref{defz} for $u$
\begin{equation}\label{ufromz}
u=\tau\mathfrak{e}_\tau\cdot u(0)+\mathfrak{e}_\tau*z 
\text{ with } \mathfrak{e}_\tau(t)=\frac{1}{\tau}e^{-\frac{t}{\tau}},
\end{equation} 
satisfies
\begin{equation}\label{wave_z_memory}
z_{tt}-(c^2+\tfrac{\delta}{\tau})\Delta z +\tfrac{\delta}{\tau} \mathfrak{e}_\tau*\Delta z
= f(t)-\delta\mathfrak{e}_\tau\cdot \Delta u(0),
\end{equation}
cf.~\cite{bucci2019regularity}.
We here also refer to \cite{dell2017moore} in which a similar relation is derived, but for the solution $u$ itself, that can be characterized by a linear viscoelastic model with an exponential kernel.


When considering the \SMGT\, equation itself with memory, one can expect some dissipation to be induced by the memory term under appropriate sign conditions. However, the degree of singularity clearly has a substantial influence on the decay behaviour. In \cite{LasieckaWang15a} 
Lasiecka and Wang consider the equation
\begin{equation}\label{MGT_memory}
\tau u_{ttt}+u_{tt}-c^2 \Delta u- b \Delta u_t -\mathfrak{k}*\Delta w = 0
\end{equation}
with three types of memory 
\begin{equation}\label{memorytypes}
w=\begin{cases}
u&\text{ type I}\\
u_t&\text{ type II}\\
\tau u_t + u&\text{ type III}
\end{cases}
\end{equation}
under the conditions
\[
\mathfrak{k}\in C^1(0,\infty)\cap C[0,\infty), \quad \mathfrak{k}\geq0,\ \mathfrak{k}'\leq0, \quad \mathfrak{k}'\leq c_0 \mathfrak{k}
\]
for some positive constant $c_0$.
For type III memory, exponential decay of the energy
\[
\begin{aligned}
\mathcal{E}_{\mathfrak{k}}[u](t):=&\nLtwo{u_{tt}}^2+\nLtwo{\nabla u_t}^2+\nLtwo{\nabla u}^2\\
&+\int_0^t\mathfrak{k}(t-s)\nLtwo{\nabla w(t)-\nabla w(s)}^2\, ds
\end{aligned}
\]
is established even in the critical case $\delta=b-\tau c^2=0$, provided $\int_0^\infty \mathfrak{k}(t)\, dt <b$. 
whereas type I memory 
requires dissipation $\delta=b-\tau c^2>0$ for leading to exponential decay.
This result is further developed in \cite{dell2016moore,LasieckaWang15b}, where general decay is studied and it is shown that with type I memory, exponential decay of the energy $\mathcal{E}_{\mathfrak{k}}[u](t)$ can only occur if the Laplacian is replaced by a bounded operator. Still, the modified energy 
\[
\begin{aligned}
\mathcal{E}_{-\mathfrak{k}'}[u](t):=&\nLtwo{u_{tt}}^2+\nLtwo{\nabla u_t}^2+\nLtwo{\nabla u}^2\\
&-\int_0^t\mathfrak{k}'(t-s)\nLtwo{\nabla u(t)-\nabla u(s)}^2\, ds
\end{aligned}
\]
must tend to zero 
as $t\to\infty$, cf. \cite{dell2016moore}.

For further results on the \SMGT\, equation with memory, we refer to, e.g., \cite{zbMATH06661873,lasiecka2017global,alves2018moore,dell2020note,Chen:2024} and the references provided therein.

The nonlinear JMGT case has first been studied in a series of papers by Nikoli\'{c} and Said-Houari: 
In \cite{NikolicSaidHouari:2021_memory_unboundeddomain}, 
\revision{local-in-time} 
\MarginRone{17.}  
 \revision{wellposedness} with possibly large initial and 
\revision{global-in-time} 
\MarginRone{22.}  
 \revision{wellposedness} as well as polynomial decay with small initial data of \eqref{JMGT-Westervelt} with type I memory on $\mathbb{R}^3$ is shown; 
\cite{NikolicSaidHouari:2021_hereditary,NikolicSaidHouari:2021_inviscid} 
on $\mathbb{R}^d$, $d\geq3$, 
extend the analysis to even including the gradient nonlinearity \eqref{JMGT-Kuznetsov} where the latter focuses on the critical case with type I memory, 
tackling the fact that the typical linear decay estimates are of regularity loss type by means of well-constructed time weigths. 

\medskip

Closely related, the role of fractional attenuation, required to model fractional power frequency dependence of damping as typical for ultrasound propagation, see, e.g., \cite{ChenHolm:2004,CaiChenFangHolm_survey2018,Szabo:1994,TreebyCox:2010,Wismer:2006}, has recently been analyzed in the context of the JMGT equation in, e.g. \cite{fracJMGT,MelianiSaidHouari:2025,Nikolic:2024_fractional}. The models under consideration arise from a substitution of the Maxwell-Cattaneo heat flux law by time fractional versions
cf.~\cite{compte1997generalized} of the form:
\begin{equation} \label{fractional_law}
	\begin{aligned}
		(1+\tau^{\alpha_1} \partial_t^{\alpha_1})\boldsymbol{q}(t) = -\kappa \partial_t^{\alpha_2} \nabla \theta,
	\end{aligned}
\end{equation} 
with the 
\revision{Djrbashian}
\MarginRone{26.}
-Caputo derivative defined by
\[
\partial_t^\alpha v = \mathfrak{k}_\alpha*v_t, \quad \mathfrak{k}_\alpha(t)=\frac{1}{\Gamma(1-\alpha)}t^{-\alpha}.
\]
These extensions are additionally motivated by the fact that the model relying on \eqref{MC_law} may violate the second law of thermodynamics~\cite{zhang2014time, fabrizio2017modeling, ferrillo2018comparing} and fractional generalizations of the heat flux law have emerged in the literature as a way of interpolating between the properties of the two flux laws \eqref{fourier_law} and \eqref{MC_law}; see, e.g.,~\cite{povstenko2011fractional, compte1997generalized, fabrizio2015some, atanackovic2012cattaneo} and the references contained therein. 
Alternatively or additionally to that, 
\revision{analogues} 
\MarginRone{27.}
of models for viscoelasticity can lead to fractional time derivatives in the model, see e.g. the fractional Zener model in \cite[Chapter 7]{frac_book}.
\subsection{Periodic solutions}\label{sec:periodic}
As 
\revision{time-periodic} 
\MarginRone{28.}
(continuous wave CW) excitations are often used in ultrasonics, the question of existence of solutions to the JMGT equation \eqref{JMGT-Westervelt} or \eqref{JMGT-Kuznetsov} under periodicity rather than initial conditions is practically relevant.
Indeed, imposing 
\[
u(T)=u(0), \quad u_t(T)=u_t(0), \quad u_{tt}(T)=u_{tt}(0)
\]
in place of initial conditions, it is immediate that when evaluating the energy identity \eqref{enid} at $t=T$, all $\Big|_0^T$ terms vanish and under the conditions \eqref{condcoeff_local}, \eqref{condcoeff_global} we end up with an estimate of the form 
\begin{equation}\label{enid_periodic}
\begin{aligned}   
&\tau\|\nabla u_{tt}\|_{L^2(0,T;L^2(\Omega))}^2+\|u\|_{H^1(0,T;H^2(\Omega))}^2\\
&\leq C\left|\int_0^t\prodLtwo{f}{-\Delta(\tau u_{tt}+\sigma u_t+\rho u)}\right|.
\end{aligned}
\end{equation} 
Using an alternative bound of the right hand side 
\begin{equation}\label{boundf_alt}
\begin{aligned}
&\int_0^t\prodLtwo{f}{-\Delta(\tau u_{tt}+\sigma u_t+\rho u)}\, ds\\
&\leq \frac{1}{2\epsilon}\int_0^t\Bigl(\tau \nLtwo{\nabla f}^2+\nLtwo{f}^2\Bigr)\, ds\\
&\quad+ \frac{\epsilon}{2}\int_0^t\Bigl(\tau \nLtwo{\nabla u_{tt}}^2+\nLtwo{\sigma\Delta u_t+\rho\Delta u}^2\Bigr)\, ds   
\end{aligned}
\end{equation}
we obtain an energy estimate of the form
\begin{equation}\label{enest_periodic_lin}
\begin{aligned}   
&\tau\|\nabla u_{tt}\|_{L^2(0,T;L^2(\Omega))}^2+\|u\|_{H^1(0,T;H^2(\Omega))}^2\\
&\leq\,C\left(\| f\|^2_{L^2(0,T;L^2(\Omega))}+\tau\|\nabla f\|^2_{L^2(0,T;L^2(\Omega))}\right). 
\end{aligned}
\end{equation}
The 
\revision{advantage} 
\MarginRone{29.}
of not needing to estimate certain terms as compared to Section~\ref{sec:ivp}  has to be paid for by the loss of $L^\infty$ in time estimates.
\revision{Using the alternative right hand side bound \eqref{boundf_alt} allows us to cope with this loss.}
\MarginRtwo{6.}
With $f=-(\eta u^2+r)_{tt}$ for a function $r$ modelling excitation and a 
\revision{fixed-point}
\MarginRone{14.} 
argument this leads to an energy estimate 
\begin{equation}\label{enest_periodic}
\begin{aligned}   
&\tau\|\nabla u_{tt}\|_{L^2(0,T;L^2(\Omega))}^2+\|u\|_{H^1(0,T;H^2(\Omega))}^2\\
&\leq\,C\left(\| r_{tt}\|^2_{L^2(0,T;L^2(\Omega))}+\tau\|\nabla r_{tt}\|^2_{L^2(0,T;L^2(\Omega))}\right). 
\end{aligned}
\end{equation}
for solutions of the nonlinear periodic problem
\begin{equation}\label{JMGT_periodic}
	\begin{aligned}
		&\tau u_{ttt}+u_{tt}-c^2 \Delta u- b \Delta u_t = -\eta(u^2)_{tt}-r_{tt} &&\text{ in }(0,T)\times\Omega\\
		&u=0 &&\text{ on }(0,T)\times\partial\Omega\\
	    &u(0)=u(T), \quad u_t(0)=u_t(T), \quad u_{tt}(0)=u_{tt}(T).		
	\end{aligned}
\end{equation}
This energy estimate, together with a Galerkin 
\revision{semidiscretization} 
\MarginRone{30.}
in space and the use of Bloch-Floquet theory for the resulting system of ODEs allows to prove \revision{wellposedness} of \eqref{JMGT_periodic} and, using a higher order energy 
\revision{estimate}, 
\MarginRone{31.}
also the periodic counterpart of \eqref{JMGT-Kuznetsov} containing gradient nonlinearities, cf. \cite{periodicJMGT}.

Another motivation for considering periodic solution is the possibility of expanding 
both the excitation and the solution of \eqref{JMGT_periodic} in multiharmonic series
\[
r(x,t)=\Re\left(\sum_{k=1}^\infty \hat{r}_k(x) e^{\imath k \omega t}\right), \quad 
u(x,t)= \Re\left(\sum_{k=1}^\infty \hat{u}_k(x) e^{\imath k \omega t}\right),
\]
(with $\omega=\frac{2\pi}{T}$), relying on existence of a $T$-periodic solution and on completeness of the functions $\cos(\omega t)$, $\sin(\omega t)$ in $L^2(0,T)$.
Inserting into \eqref{JMGT_periodic} leads to a discrete convolution resulting from the squared terms and allows us to equivalently formulate the equation in frequency domain as a coupled system of Helmholtz type equations
\begin{equation}\label{multiharmonic}
\begin{aligned}
&[\imath m^3\omega^3\tau+m^2\omega^2
+ c^2\Delta +(\tau c^2+\delta) \imath m\omega \Delta]\hat{u}_m\\
&=\frac12 \, \eta\, m^2\omega^2\, \left(\sum_{\ell=1}^{m-1} \hat{u}_\ell \hat{u}_{m-\ell}
+ \sum_{k=1}^\infty\overline{\hat{u}_{k}} \hat{u}_{k+m} 
+ \sum_{k=1}^\infty\hat{u}_{m+k} \overline{\hat{u}_{k}}\right)
+ m^2\omega^2\, \hat{r}_m\\
&m\in\mathbb{N}.
\end{aligned}
\end{equation}
This can also serve as a mathematical explanation of the appearance of so-called higher harmonics, that is, contributions at multiples of the fundamental frequency $\omega$ even if the excitation is only applied at frequency $\omega$, that is, $\hat{r}_m=0$, $m\geq 2$.
  
\subsection{Relation to other models of nonlinear acoustics via singular limits}\label{sec:tau2zero}
To rigorously connect \eqref{JMGT-Westervelt} and \eqref{JMGT-Kuznetsov} to the classical  (Lighthill-)Westervelt~\cite{lighthill1956viscosity,westervelt1963parametric} and Kuznetsov~\cite{kuznetsov1971equations} models, respectively, as their $\tau=0$ limiting cases, in \cite{bongarti2020vanishing,JMGT_BKVN,JMGT_Neumann} an analysis of convergence in appropriate function spaces has been carried out.
The challenges in such an analysis arise from the fact that as $\tau\to0$ the PDE changes its type from hyperbolic to parabolic: As mentioned above, the first order formulation of its linearization in case $\tau>0$ gives rise to a group, whereas the $\tau=0$ case is known to lead to an analytic semigroup and maximal parabolic regularity.

An essential ingredient for such a limiting 
\revision{analysis} 
\MarginRone{32.}
are $\tau$ uniform bounds.
To this end, in the derivation of \eqref{energy_est}, exploiting dissipativity $\delta>0$, we employ an alternative estimate of the right hand side 
\begin{equation*}
\begin{aligned}
&\int_0^t\prodLtwo{f}{-\Delta(\tau u_{tt}+\sigma u_t+\rho u)}\, ds\\
&\leq \frac{1}{2\epsilon}\int_0^t\Bigl(\tau \nLtwo{\nabla f}^2+\nLtwo{f}^2\Bigr)\, ds\\
&\quad+ \frac{\epsilon}{2}\int_0^t\Bigl(\tau \nLtwo{\nabla u_{tt}}^2+\nLtwo{\sigma\Delta u_t+\rho\Delta u}^2\Bigr)\, ds.   
\end{aligned}
\end{equation*}
Absorbing the $\epsilon$ terms in 
\revision{dissiative} 
\MarginRone{33.}
terms on the \revision{left-hand-side} \MarginRone{19.} of \eqref{enid} and bootstrapping the $L^2(0,t;L^2(\Omega))$ norm of the third order in time term from the PDE thus yields an energy estimate of the form
	\begin{equation*}
	\begin{aligned}
\|u^\tau\|_{W,\tau}^2:=
& \revision{
\tau^2 \|u^\tau_{tt}\|^2_{L^\infty(0,T;H^1(\Omega))}+\tau \|u^\tau_t\|^2_{L^\infty(0,T;H^2(\Omega))}+\|u^\tau\|^2_{L^\infty(0,T;H^2(\Omega))}
}\\
&  +\tau^2 \|u^\tau_{ttt}\|^2_{L^2(0,T;L^2(\Omega))} + \tau \|u^\tau_{tt}\|^2_{L^2(0,T;H^1(\Omega))}+\|u^\tau\|^2_{H^1(0,T;H^2(\Omega))}
\\
\leq&\,C\Bigl(\revisionown{\tau^2|u_2|^2_{H^1(\Omega)}+\tau|u_1|^2_{H^2(\Omega)}+|u_0|^2_{H^2(\Omega)}}\\
&\qquad+ \| f\|^2_{L^2(0,T;L^2(\Omega))}+\tau\|\nabla f\|^2_{L^2(0,T;L^2(\Omega))}\Bigr). 
	\end{aligned}
	\end{equation*}
\MarginRone{f}
with a constant $C$ independent of $T$ and $\tau$ for solutions of the linearized equation \eqref{JMGT_lin}. 
By a 
\revision{fixed-point}
\MarginRone{14.} 
argument and with small initial data this can be carried over to the nonlinear setting \eqref{JMGT-Westervelt} as
	\begin{equation*}
	\begin{aligned}
\|u^\tau\|_{W,\tau}^2\leq C\left(|u_0|^2_{H^2(\Omega)}+\tau|u_1|^2_{H^2(\Omega)}+\tau|u_2|^2_{H^1}\right)
	\end{aligned}
	\end{equation*}
Thus, uniform boundedness of the $\tau$ independent part $\|u^\tau\|^2_{H^1(0,T;H^2(\Omega))}$ (and, again, by a bootstrapping argument, also of $\|u^\tau_{tt}\|^2_{L^2(0,T;L^2(\Omega))}$) can be established and yields weak convergence in $H^2(0,T;L^2(\Omega))\cap H^1(0,T;H^2(\Omega))$ to an element $\bar{u}$ in this space, which can be shown to solve the $\tau=0$ equation.

Note that we have here employed a function space setting and energy estimates that are slightly different from \cite{bongarti2020vanishing,JMGT_BKVN,JMGT_Neumann}, in order to allow for a unified exposition re-using \eqref{enid}.

Similar studies have been carried out including the quadratic gradient nonlinearity \eqref{JMGT-Kuznetsov}, as well as for fractional attenuation and in the 
\revision{time-periodic} 
\MarginRone{28.}
setting, cf.~\cite{frac_tau2zero_PartII,Nikolic:2024_fractional,periodicJMGT}. 

On the other hand, the inviscid limit $\delta\searrow0$ towards the critical case for fixed positive $\tau>0$ has been analyzed in \cite{b2zeroJMGT}. 
\section{Some control and inverse problems}\label{sec:controlinverse}
\subsection{Control}\label{subsec:control}
\revision{In \cite{bucci2018feedback},} 
\MarginRone{34.}
Bucci and Lasiecka study the problem of controlling the acoustic pressure $u$ governed by the \SMGT\, equation~\eqref{MGT} such that a desired pressure distribution $u^d$ is followed as tightly as possible.
This is formulated as the minimization of an objective function consisting of a tracking term and a control cost 
\begin{equation}\label{J}
J(g,u)=\int_0^T\int_\Omega|u-u^d|^2\, dx\, dt +\gamma \int_0^T\int_{\Gamma_0} |g|^2\, dS\, dt.
\end{equation}
with some positive cost parameter $\gamma>0$.
As relevant for practical applications of e.g., high intensity ultrasound, the control function $g$ acts on a part of  the boundary $\Gamma_0\subseteq\partial\Omega$ in the form of a Neumann boundary condition, modelling excitation by, e.g., an array of piezoelectric transducers. The rest of the boundary is equipped with absorbing boundary conditions to model nonreflecting boundary conditions, 
\revision{i.e., free propagation of waves through} 
\MarginRone{35.}
that boundary part.
Note that in \eqref{J}, only the $L^2$ with respect to time norm of the control is used, in order to allow for non-smooth (e.g., switching) controls.
The problem of minimizing $J$ subject to the PDE constraint 
\begin{equation}\label{MGT-Neumann}
	\begin{aligned}
		&\tau u_{ttt}+u_{tt}-c^2 \Delta u- b \Delta u_t = 0 &&\text{ in }(0,T)\times\Omega\\
		&\partial_\nu u=g &&\text{ on }(0,T)\times\Gamma_0\\
		&c\partial_\nu u+u_t=0 &&\text{ on }(0,T)\times\Gamma_1=\partial\Omega\setminus\Gamma_0\\
	    &u(0)=u_0, \quad u_t(0)=u_1, \quad u_{tt}(0)=u_2		
	\end{aligned}
\end{equation}
\MarginRone{36.}
comes with several challenges.
Firstly, the appearance of unbounded control operators due to the use of a boundary control is here not counteracted by a regularizing effect of the evolution dynamics; this is due to the hyperbolic nature of the S/J MGT equation as compared to the parabolic one of strongly damped second order wave equations (such as the $\tau=0$ case of the classical Westervelt and Kuznetsov equations).
Secondly, as a consequence of the $b$ term, extension of the Neumann data into the interior of $\Omega$ also involves the time derivative of the control, but the objective \eqref{J} obviously lacks coercivity with respect to $g_t$.
To cope with this, a control-to-state map relying on the variation of parameters formula for the semigroup governing the first order reformulation of \eqref{MGT} (see also \cite{kaltenbacher2011wellposedness,marchand2012abstract}) is derived and substantiated in \cite{bucci2018feedback}, that relies on values of $g(t)$, but not on $g_t(t)$.
By adding a correction term involving the initial values of the control, 
\revision{wellposedness} of the minimization problem is achieved. 
Moreover, the optimal control $g$ is given in feedback form, that is, the control at each time instance is expressed in terms of the state at that time instance, via some linear 
\revision{time-dependent} 
\MarginRone{37.}
feedback operator (that is, a so-called feedback synthesis is established). The 
\revision{time-dependent} 
\MarginRone{37.}
part of this feedback operator is 
\revision{governed}  
\MarginRone{38.}
by a non-standard Riccati equation, that is also shown to be well-posed in \cite{bucci2018feedback}.
\subsection{Imaging with nonlinear ultrasound}\label{subsec:nonlinearityimaging}
\revision{Model-based} 
\MarginRone{39.}
quantitative tomography relies on the fact that certain coefficients contained in the PDE model are specific to the tissue type. Therefore, maps of these coefficients as functions of the spatial variables provide medical imaging tools and in fact contain clinically useful information beyond these images.
In the JMGT equation as a model of nonlinear ultrasound, the relevant quantities are the sound speed $c$, the attenuation $b$, and the nonlinearity coefficient $\eta$.
The above mentioned imaging task amounts to reconstructing these coefficients as functions of space from additional observations -- typically measurements of the pressure $u$ 
\begin{equation}\label{observation}
p^{obs}(t,x_0) = u(t,x_0), \quad(t,x_0)\in(0,T)\times \Sigma
\end{equation}
where $\Sigma\subseteq\overline{\Omega}$ is a smooth $d-1$ dimensional manifold, modeling, e.g. an array of piezoelectric transducers or hydrophones.
A crucial question for this inverse problem is uniqueness, that is, whether the given data \eqref{observation} suffices to uniquely determine the 
\revision{sought-after} 
\MarginRone{40.}
quantities as functions of the $d$ space variables in $\Omega$.

Here it turns out that nonlinearity helps: The fact that even when excited at a single frequency $\omega$, the response contains contributions at all multiples of $\omega$ according to \eqref{multiharmonic}, illustrates the multiplication of information due to nonlinearity.
An additional positive effect of the relaxation time term lies in its role of re-establishing a hyperbolic character of the PDE and finite speed of propagation, thus counteracting the loss of information caused by strong attenuation.
As a consequence, besides the JMGT equation being a better physical model, it also allows to mathematically prove better uniqueness results.
Indeed, it has been shown in \cite{nonlinearity_imaging_JMGT} for the initial value problem and in \cite{nonlinear_imaging_JMGT_freq} for the 
\revision{time-periodic} 
\MarginRone{28.}
setting, that observation of $u$ over time \eqref{observation} from a single source suffices for local uniqueness of the nonlinearity coefficient $\eta(x)$. Adding a second observation by just modifying the amplitude of the first source allows one to recover the sound speed $c(x)$ as well \cite{nonlinearity_imaging_JMGTmulticoeff}. This is in stark contrast to the linear setting, where it is known that infinitely many sources and observations (the full Dirichlet-Neumann map for the underlying PDE) are needed for guaranteeing uniqueness of $c(x)$.

To give an idea of such a uniqueness proof and the role of the relaxation term therein, we sketch the arguments from \cite{nonlinearity_imaging_JMGTmulticoeff}, but restrict exposition  to the reconstruction of $\eta(x)$ alone, for the sake of transparency.
To do so, we 
\def\trSigma{\text{tr}_{\Sigma}}
\newcommand{\uld}[1]{\underline{d{#1}}}
write the inverse problem as an operator equation for the unknown functions $\eta=\eta(x)$ and $u=u(t,x)$ 
\[
F(\eta,u)=y
\]
with 
\[
F(\eta,u)=\left(\begin{array}{l}
\mathcal{L} u + \eta u^2\\
\trSigma u \end{array}\right), \quad 
y=\left(\begin{array}{l}
r^{src}\\
p^{obs}\end{array}\right)
\]
with the linear differential-integral operator defined by 
\[
\mathcal{L}u:=\tau u_t+u-c^2\Delta \int_0^t\int_0^s u(r)\, dr\, ds -b \Delta \int_0^t u(s) \, ds,
\]
a given source of the form $f=-r^{scr}_{tt}$ cf. \eqref{JMGT-Westervelt}, and given observations $p^{obs}$ cf. \eqref{observation}.  
Time periodicity conditions $u(T)=u(0)$ are understood to be incorporated into the function space used for $u$ and boundary conditions (e.g., Dirichlet ones) are included in the definition of $-\Delta$.

We linearize $F$ at a reference point $(\eta^0,u^0)$ that we choose as $\eta^0=0$ and $u^0(t,x)$ of space-time separable form $u^0(t,x)=\psi(t)\phi(x)$, so that the linearized equation becomes
\[
F'(0,u^0)(\uld{\eta},\uld{u})=\left(\begin{array}{l}
\mathcal{L} \uld{u} + \psi^2\,\phi^2\,\uld{\eta}\\
\text{tr}_{\Sigma} \uld{u} \end{array}\right)
=\left(\begin{array}{l}
\uld{r}\\
\uld{p}\end{array}\right)
\]
Re-defining $\uld{\widetilde{\eta}}(x)=\phi(x)^2\uld{\eta}(x)$, we can formally resolve this equation as 
\[
\begin{aligned}
&\uld{u}= \mathcal{L}^{-1}\left(\uld{r}-\psi^2\uld{\widetilde{\eta}}\right)\\
&\uld{\widetilde{\eta}}=\psi^{-2}\cdot\Bigl(\text{tr}_{\Sigma} \mathcal{L}^{-1}\Bigr)^{-1}\Bigl(\text{tr}_{\Sigma}\mathcal{L}^{-1}\uld{r}-\uld{p}\Bigr)\end{aligned}
\]
In order to substantiate this formula, it is important to keep in mind that the trace operator is clearly not injective -- however its restriction to an eigenspace of the negative Dirichlet Laplacian $-\Delta$ is, by unique continuation under quite general assumptions on the observation surface $\Sigma$ see, e.g. \cite{
JiangLiPauronYamamoto2023,
Tolsa:2023} and the references therein. 
To disentangle the eigenspaces for this purpose, we rely on the frequency domain formulation induced by \eqref{multiharmonic}. 
\revision{We analytically}  
\MarginRone{41.}
extend the observation in frequency domain from the discrete set $\{\imath m \omega\, : \, m\in\mathbb{N}\}$ to the whole complex plane, up to a countable set of poles that are the zeros of the function $z\mapsto\tau z^3+z^2+(c^2+bz)\lambda_j$ with $\lambda_j$ eigenvalue of $-\Delta$. Considering the residues at these poles allows us to single out the eigenspace contribution corresponding to $\lambda_j$.  
This allows one to make sense of $\text{tr}_{\Sigma}\mathcal{L}^{-1}$ and its inverse in the formula above.

With a natural choice of the topology in preimage space, e.g. 
$(\uld{\widetilde{\eta}},\uld{u})\in H^1(\Omega)\times L^2(0,T;H^1(\Omega))$ and defining the topology in image space by 
$\|(\uld{r},\uld{p})\|_Y:=$\\ $\|F'(0,u^0)^{-1}(\uld{r},\uld{p})\|_{H^1(\Omega)\times L^2(0,T;H^1(\Omega))}$ trivially renders $F'(0,u^0)$ an isomorphism and the Inverse Function Theorem yields local uniqueness and even stability. A refined analysis establishes bounds of this artificially defined data space norm by means of Sobolev norms, with a constant that depends on $\tau>0$ and that blows up as $\tau\searrow0$. 
This shows the importance of the JMGT model (as compared to the classical $\tau=0$ one) also for this inverse problem application.

\MarginRtwo{7.}
\revision{
\section{Some open problems} 
We conclude by listing a small (and again, incomplete) selection of open problems.
\begin{itemize}
\item 
While the {\em critical (inviscid) case $\delta=0$} is already well understood for initial value problems for the JMGT equation, it appears to be still open in the {\em time periodic} setting. 
\item
Instead of combining the fundamental balance and constitutive equations to a single third order in time equation, one could also consider an analysis of the {\em first order system} resulting from these equations, as has been done in, e.g.,  \cite{DekkersRozanovaKhodygo2020,QuadraticWave} for other models of nonlinear acoustics.
\item
The formal calculations in Section~\ref{subsec:energyestimates} show that the auxiliary state $z$ up to a lower order remainder term satisfies a weakly damped wave equation, and thus clearly provide evidence for the physical justification of {\em finite speed of propagation}.  
A mathematically rigorous verification of the resulting conjecture on a (generalized) cone of dependence, that is, the fact that for each positive time instance, the subset on which the initial data uniquely determine the JMGT solution is bounded, is still to be done. 
\end{itemize}
}
\section*{Acknowledgment}
This research was funded in part by the Austrian Science Fund (FWF) 
[10.55776/P36318]. 
For open access purposes, the author has applied a CC BY public copyright license to any author accepted manuscript version arising from this submission.
%
\MarginRone{42.}
\end{document}